\documentclass[
10.5pt]{article}
\usepackage{amsmath,amssymb,color}
\def\R{{\mathbb{R}}}
\def\Z{{\mathbb{Z}}}
\def\BMO{{\mathrm{BMO}}}

\def\Lip{{\mathrm{Lip}}}

\def\Osc{{\mathrm{Osc}}}
\def\dint{\displaystyle\int}
\def\dsup{\displaystyle\sup}

\usepackage{graphicx,amssymb,mathrsfs,amsmath,latexsym,amsfonts,titlesec,amscd,epsfig,cite,amsthm}
\usepackage{indentfirst} 

\DeclareMathOperator*{\supp}{supp}

\newcommand{\upcite}[1]{\textsuperscript{\textsuperscript{\cite{#1}}}}

\theoremstyle{definition} 

\newtheorem{theorem}{\indent
                  Theorem}[section]

    \newtheorem{lemma}{\indent  Lemma} [section]

\theoremstyle{definition} 
        \newtheorem{remark}{\indent  Remark}  

\newtheorem{theoremalph}{\indent Theorem}
\title{\bf\Large
Weighted Endpoint Estimates for Multilinear Commutators of Marcinkiewicz Integrals }
\author{
{\normalsize Jianglong WU$^{1}$ \  Qingguo LIU$^{2}$ 
} 
\\
(\small 1 College of Science, Mudanjiang Normal University, Mudanjiang, 157012, China)\\
(\small 2 Laboratory for Multiphase Processes, University of Nova Gorica, Nova Gorica, SI-5000 , Slovenia)\\
 \small 
}
\date{} 
\begin{document}

\maketitle

\noindent

\begin{minipage}[t]{15cm}

\setlength{\baselineskip}{1.0em}

\noindent  

 { \bf Abstract:} Let $\mu_{\Omega,\vec{b}}$ be the
multilinear commutator generalized by $\mu_{\Omega}$, the
$n$-dimensional Marcinkiewicz integral, with $\Osc_{\exp
L^{^{\tau}}}(\R^{n})$ functions for $\tau\ge 1$, where $\Osc_{\exp
L^{^{\tau}}}(\R^{n})$
is a space of Orlicz type satisfying that $\Osc_{\exp L^{^{\tau}}}(\R^{n})=\BMO(\R^{n})$ if $\tau=1$ and $\Osc_{\exp L^{^{\tau}}}(\R^{n})\subset\BMO(\R^{n})$ if $\tau>1$. 
The authors establish the weighted weak $L\log L$-type estimates for
$\mu_{\Omega,\vec{b}}$ when $\Omega$ satisfies a kind of Dini
conditions.
\\

 \\

 { \bf MR(2000) Subject Classification:}  \ 42B20; 42B25; 42B99  \ \

\end{minipage}


\section{Introduction and Main Result}

Denote by $S^{n-1}$ the unit sphere in $\R^{n}~(n\ge 2)$ equipped
with the normalized Lebesgue measure
$\mathrm{d}x'=\mathrm{d}\sigma(x')$. Let $\Omega(x)\in
L^{1}(S^{n-1})$ be homogeneous function of degree zero in $\R^{n}$
satisfying
\begin{equation} \label{equ:1}
\dint_{S^{n-1}} \Omega(x') \mathrm{d}x'=0,
\end{equation}
where $x'=x/|x|~(x\neq 0)$.

The n-dimensional Marcinkiewicz integral introduced  by Stein
\upcite{S1} is defined by
$$\mu_{\Omega}(f)(x)=\Big(\dint_{0}^{\infty} \Big| \dint_{|x-y|\le t} \frac{\Omega(x-y)}{|x-y|^{n-1}} f(y) \mathrm{d}y\Big|^{2} \frac{\mathrm{d}t}{t^{3}}\Big)^{\frac{1}{2}},\ \ \ x\in \R^{n}.$$

A weight will always means a positive locally
 integrable function. As usual, we denote by $A_{p}~(1\le p \le
 \infty)$ the Muckenhoupt weights classes (see \cite{S2} and \cite{G} for
 details). For a weight $\omega$ on $\R^{n}$, we write $\|f\|_{L^{p}_{\omega}(\R^{n})}= (\int_{\R^{n}} |f(x)|^{p} \omega(x) \mathrm{d}x)^{1/p}$
 and $\omega(E)=\int_{E}\omega(x) \mathrm{d}x$.

 In 2004, Ding, Lu and Zhang\upcite{DLZ} studied the weighted weak
$L\log L$-type estimates for the commutators of the Marcinkiewicz
integral, which is defined by
$$\mu_{\Omega,b}^{m}(f)(x)=\Big(\dint_{0}^{\infty} \Big| \dint_{|x-y|\le t} \frac{(b(x)-b(y))^{m}\Omega(x-y)}{|x-y|^{n-1}} f(y) \mathrm{d}y\Big|^{2} \frac{\mathrm{d}t}{t^{3}}\Big)^{\frac{1}{2}},\ \ \ m\in \Z^{+},b\in\BMO(\R^{n}),$$
when the kernel $\Omega$ satisfies the $\Lip_{\alpha} (0<\alpha\le
1)$ condition, that is, there exists a constant $C>0$ such that
\begin{equation} \label{equ:02}
 |\Omega(x')-\Omega(y')| \le C|x'-y'|^{\alpha}, \ \ \forall\ x', y'\in S^{n-1}.
 \end{equation}

In 2008, Zhang\upcite{Z1} established the weighted weak $L(\log
L)^{1/r}$-type estimates for the multilinear commutators of the
Marcinkiewicz integral when $\omega\in A_{1}$, and $\Omega$
satisfies \noindent(\ref{equ:1}) and \noindent(\ref{equ:02}).


Let $\Omega\in L^{r}(S^{n-1})~(r\ge 1)$, the integral modulus of
continuity of order $r$ of $\Omega$ is defined by
$$\omega_{r}(\delta)=\dsup_{|\rho|<\delta} \Big( \dint_{S^{n-1}} |\Omega(\rho x')-\Omega(x')|^{r} \mathrm{d}x'\Big)^{1/r},$$
where $\rho$ is a rotation in $\R^{n}$ with $|\rho|=\sup_{x'\in
S^{n-1}} |\rho x'-x'|$.

We say $\Omega\in L^{r}(S^{n-1})~(r\ge 1)$ satisfies the
$L^{r}$-Dini condition if $\int_{0}^{1} \omega_{r}(\delta)
\delta^{-1} \mathrm{d}\delta<\infty$.

Recently, Zhang\upcite{Z2} also considered  the following result. 

\begin{theoremalph}\upcite{Z2}\label{thm.A}
Let $b\in \BMO(\R^n), \Omega\in L^{r}(S^{n-1})$ for some $r>1$, and
$\omega^{r'}\in A_{1}$. If $\Omega$ satisfies \noindent(\ref{equ:1})
and
\begin{equation} \label{equ:2}
\dint_{0}^{1} \frac{\omega_{r}(\delta)}{\delta} \Big(\log
\frac{1}{\delta}\Big)^{m}
 \mathrm{d}\delta<\infty,
 \end{equation}
then for all $\lambda>0$, there has
$$ \omega(\{x\in \R^{n}: {\color{red} \mu_{\Omega,b}^{m} }(f)(x)> \lambda\}) \le C  
\dint_{\R^{n}} \frac{|f(y)|}{\lambda}
\Big(1+\log^{+}\frac{|f(y)|}{\lambda} \Big)^{m} \omega(y)
\mathrm{d}y,$$ where $C$ is a positive constant independent of $f$
and $\lambda$.
\end{theoremalph}

In this paper, by applying the calder\'{o}n-Zygmund decomposition
theory, we will study the weighted weak $L\log L$-type estimates for
the multilinear commutators generated by $\mu_{\Omega}$  and
$\Osc_{\exp L^{^{r}}}(\R^{n})$ functions, in analogy with the
results established by P\'{e}rez and Trujillo-Gonz\'{a}lez in
\cite{PT} for the multilinear commutators of Calder\'{o}n-Zygmund
operators. Before stating our results, we first recall some
notation.

 Let $ m$ be a positive
integer and $\vec{b}=(b_{1},b_{2},\cdots,b_{m})$, we define the
multilinear commutators $\mu_{\Omega,\vec{b}}$   by
$$\mu_{\Omega,\vec{b}}(f)(x)=\Big(\dint_{0}^{\infty} \Big| \dint_{|x-y|\le t} \frac{\Omega(x-y)f(y) }{|x-y|^{n-1}} \prod_{j=1}^{m}\Big(b_{j}(x)-b_{j}(y)
\Big) \mathrm{d}y\Big|^{2}
\frac{\mathrm{d}t}{t^{3}}\Big)^{\frac{1}{2}}.$$ It is easy to see,
when $m = 1$, $\mu_{\Omega,\vec{b}}$ is the commutator of
Marcinkiewicz integral and when $b_{1}=\cdots=b_{m}$,
$\mu_{\Omega,\vec{b}}$ is the higher commutator of Marcinkiewicz
integral.

To state the weak type estimate for the multilinear commutator
$\mu_{\Omega,\vec{b}}$, we need to introduce the following notation.
As in \cite{PT}, given any positive integer $m$, for all $1\le j \le
m$, we denote by $\mathscr{C}_{j}^{m}$ the family of all finite
subsets $\sigma=\{\sigma(1),\sigma(2),\cdots, \sigma(j)\}$ of
$\{1,2,\cdots,m\}$ with $j$ different elements. For any $\sigma\in
\mathscr{C}_{j}^{m}$, we define the complementary sequence
$\sigma'=\{1,2,\cdots,m\}\setminus \sigma$.

In the following, we will always assume that $\Omega$ be homogeneous
function of degree 0, and let $\vec{b}=(b_{1},b_{2},\cdots,b_{m})$
be a finite family of locally integrable functions. For all $1\le j
\le m$ and ~$\sigma=\{\sigma(1),\sigma(2),\cdots, \\
 \sigma(j)\} \in
\mathscr{C}_{j}^{m}$, we write for any $i$-tuple $
(\tau_{_{1}},\tau_{_{2}},\cdots,\tau_{_{m}})$ with
$\tau_{_{j}}\ge 1$ for $1\le j \le m$,
$1/\tau_{_{\sigma}}=1/\tau_{_{\sigma(1)}}+\cdots+1/\tau_{_{\sigma(j)}}$
and $1/\tau_{_{\sigma'}}=1/\tau-1/\tau_{_{\sigma}}$, where
$1/\tau=1/\tau_{_{1}}+\cdots+1/\tau_{_{m}}$, we will denote
$\vec{b}_{\sigma}=(b_{\sigma(1)},b_{\sigma(2)},\cdots
,b_{\sigma(j)})$ and the product
$b_{\sigma}=b_{\sigma(1)}b_{\sigma(2)}\cdots b_{\sigma(j)}.$ With
this notation, we write
$$\|\vec{b}_{\sigma}\|_{\Osc_{\exp L^{\tau_{_{\sigma}}}} (\R^{n})}= 
\|b_{\sigma(1)}\|_{\Osc_{\exp L^{\tau_{_{\sigma(1)}}}} (\R^{n})}
 \cdots \|b_{\sigma(j)}\|_{\Osc_{\exp L^{\tau_{_{\sigma(j)}}}} (\R^{n})}.$$
In particular, we write 
$$
 \big( b(x)-b(y) \big)_{\sigma}=\big( b_{\sigma(1)}(x)-b_{\sigma(1)}(y) \big) \cdots \big( b_{\sigma(j)}(x)-b_{\sigma(j)}(y)
 \big),
$$
and
$$
 \big( b_{B}-b(y) \big)_{\sigma}=\big( (b_{\sigma(1)})_{B}-b_{\sigma(1)}(y) \big) \cdots \big( (b_{\sigma(j)})_{B}-b_{\sigma(j)}(y)
 \big),
$$
where $B$ is any ball in $\R^{n}$, $x,y\in \R^{n}$,  and
$f_{B}=\frac{1}{|B|}\int_{B} f(y)\mathrm{d}y$ . For any $\sigma\in
\mathscr{C}_{j}^{m}$, we set
$$\mu_{\Omega,\vec{b}_{\sigma}}(f)(x)=\Big(\dint_{0}^{\infty} \Big| \dint_{|x-y|\le t} \frac{\Omega(x-y)f(y) }{|x-y|^{n-1}} \prod_{i=1}^{j}\Big(b_{\sigma(i)}(x)-b_{\sigma(i)}(y)
\Big) \mathrm{d}y\Big|^{2}
\frac{\mathrm{d}t}{t^{3}}\Big)^{\frac{1}{2}}.$$ If
$\sigma=\{1,2,\cdots,m\}$, then $\sigma'$ is an empty set, we
understand $\mu_{\Omega,\vec{b}_{\sigma}}=\mu_{\Omega,\vec{b}}$ and
$\mu_{\Omega,\vec{b}_{\sigma'}}=\mu_{_{\Omega}}$.

 Our result can be stated as follows.

\begin{theorem}\label{thm.1}
Let $b_{j}\in \Osc_{\exp L^{\tau_{_{j}}}}, \tau_{_{j}}\ge 1~ (1\le j
\le m), \Omega\in L^{r}(S^{n-1})$ for some $r>1$, and
$\omega^{r'}\in A_{1}$. If $\Omega$ satisfies \noindent(\ref{equ:1})
and \noindent(\ref{equ:2}),
then for all $\lambda>0$, there has
$$ \omega(\{x\in \R^{n}: \mu_{\Omega,\vec{b}}(f)(x)> \lambda\}) \le C  
\dint_{\R^{n}} \frac{|f(y)|}{\lambda}
\Big(1+\log^{+}\frac{|f(y)|}{\lambda} \Big)^{m} \omega(y)
\mathrm{d}y,$$ where $C$ is a positive constant independent of $f$
and $\lambda$.

\end{theorem}

\begin{remark}
Noting that $\Osc_{\exp L^{1}}=\BMO$ and  $\Osc_{\exp
L^{\tau}}\subset\BMO$ for $\tau>1$. For more information on Orlicz
space see \cite{RR}.
\end{remark}

Obviously, condition  \noindent(\ref{equ:2}) is slightly stronger
than the $L^{r}$-Dini condition, but much more weaker than the $\Lip
_{\alpha}$ condition. Noting that $\mu_{\Omega,\vec{b}}$ coincides
with $\mu_{\Omega,b}^{m}$ when $b_{j}=b$ for $j=1,2,\cdots,m$. So,
Theorem \ref{thm.1} improves the main results in \cite{Z1} and
\cite{Z2}.

Throughout this paper, $C$ denotes a constant that is independent of
the main parameters involved but whose value may differ from line to
line. For any index $p \in [1, \infty]$, we denote by $p'$ its
conjugate index, namely, $1/p+1/p' = 1$. For $A \sim B$, we mean
that there is a constant $C > 0$ such that $C^{-1}B \le A \le CB$.

\section{Preliminaries  and Lemmas} \setcounter{equation}{0}

In this section, we will formulate some lemmas and preliminaries.

\begin{lemma}\upcite{DL}\label{lem.1}
Suppose that $0< \alpha<n, r>1$ and $\Omega$ satisfies the
$L^{r}$-Dini condition. If there is a constant $C_{0}$ with $0<C_{0}
<1/2$ such that $|y|<C_{0}K$, then
$$ \Big( \dint_{K<|x|<2K} \Big| \frac{\Omega(x-y)}{|x-y|^{n-\alpha}}- \frac{\Omega(x)}{|x|^{n-\alpha}}\Big|^{r} \mathrm{d}x\Big)^{1/r}\le
 CK^{n/r-n+\alpha} \Big( \frac{|y|}{K} + \dint_{|y|/(2K)<\delta<|y|/K} \frac{\omega_{r}(\delta)}{\delta}\mathrm{d}\delta\Big).$$

\end{lemma}

\begin{lemma}\upcite{FS}\label{lem.2}
Suppose $\Omega\in L^{r}(S^{n-1})$ for some $r>1$ and
$\omega^{r'}\in A_{1}$. Then for any $\lambda>0$, there is a
constant $C>0$ independent of $f$ and $\lambda$, such that
$$ \omega(\{x\in \R^{n}: \mu_{\Omega}(f)(x)> \lambda\}) \le C \lambda^{-1}\|f\|_{L^{1}_{\omega}(\R^{n})}.$$
\end{lemma}

\begin{lemma}
\label{lem.3}
Let  $\omega\in A_{1}, 1< p< \infty$, $b_{j}\in \Osc_{\exp
L^{\tau_{_{j}}}},~ \tau_{_{j}}\ge 1 ~( 1\le j \le m), \Omega\in
L^{r}(S^{n-1})$ for some $r>1$ and satisfies \noindent(\ref{equ:1})
and \noindent(\ref{equ:2}). Then,  there is a constant $C>0$
independent of $f$, 
 such that
$$ \| \mu_{\Omega,\vec{b}}(f)\|_{L_{\omega}^{p}(\R^{n})}\le C \|\vec{b}\|_{\Osc_{\exp L^{\tau}}} \| f\|_{L_{\omega}^{p}(\R^{n})}.$$
\end{lemma}

The idea of  the proof of  Lemma \ref{lem.3} comes from  the
corollary 1.3 in \cite{Z1}.     We omit the details. 

We also need a few facts of Orlicz spaces, see \cite{RR} for more
information.

A function $\varphi: [0,+\infty)\to [0,+\infty)$ is called a Young
function if $\varphi$ is continuous, convex and increasing with
$\varphi(0)=0$ and $\varphi(t)\to +\infty$ as $t\to +\infty$. We
defined the $\varphi$-average of a function $f$ over a ball $B$ by
means of the Luxemburg norm
$$ \|f\|_{\varphi,B}=\inf \Big\{ \lambda>0: \frac{1}{|B|}\dint_{B} \varphi \Big(\frac{|f(y)|}{\lambda} \Big) \mathrm{d}y\le 1 \Big\},$$
which satisfies the following inequalities (see \cite{RR}, P.69 or
formula (7) in \cite{PP})
\begin{equation} \label{equ:3}
 \|f\|_{\varphi,B}\le \inf \Big\{ \eta+ \frac{\eta}{|B|}\dint_{B} \varphi \Big(\frac{|f(y)|}{\eta} \Big) \mathrm{d}y\le 1 \Big\}\le 2\|f\|_{\varphi,B}.
 \end{equation}

The Young function that we are going to be using is
$\Phi_{\alpha}(t)=t(1+\log^{+}t)^{\alpha}~ (\alpha>0)$ with its
complementary Young function $\tilde{\Phi}_{\alpha}(t)\approx
\exp(t^{1/\alpha})$. Denote by $\|f\|_{L(\log
L)^{\alpha},B}=\|f\|_{\Phi_{\alpha},B}$ and $\|f\|_{\exp
L^{1/\alpha},B}=\|f\|_{\tilde{\Phi}_{\alpha},B}$. When $\alpha=1$,
we simply denote by $\Phi(t)=t(1+\log^{+}t)$ and
$\tilde{\Phi}(t)\approx e^{t}$, and by $\|f\|_{L\log
L,B}=\|f\|_{\Phi,B}$ and $\|f\|_{\exp L,B}=\|f\|_{\tilde{\Phi},B}$.
By the generalized H\"{o}lder's inequality (see \cite{O}), we have
\begin{equation} \label{equ:4}
 \frac{1}{|B|} \dint_{B} |f(y)g(y)|\mathrm{d}y \le 2\|f\|_{L(\log L)^{\alpha},B}\|g\|_{\exp L^{1/\alpha},B}.
 \end{equation}

As usual, for a locally integrable function $f$ and a ball $B$, we
denote $f_{B}=\frac{1}{|B|}\int_{B} f(y)\mathrm{d}y$. Let $b\in
\BMO(\R^{n})$, for any ball $B$ and integer $k\ge 0$, there has (see
\cite{S2}, p.141)
\begin{equation} \label{equ:5}
|b_{2^{k+1}B}-b_{B}|\le C(k+1)\|b\|_{\ast},
 \end{equation}
where $\ell B$ denotes the $\ell$-times concentric expansion of $B$
and $\|b\|_{\ast}$ denotes the $\BMO$ norm of $b$.

By the John-Nirenberg's inequality, it is not difficult to see that
(c.f. \cite{P}, p.169)
\begin{equation} \label{equ:6}
\|b-b_{B}\|_{\exp L,B}\le C\|b\|_{\ast}.
 \end{equation}

Let $M_{L(\log L)^{\alpha}}(f)(x)=\sup_{B\ni x}\|f\|_{L(\log
L)^{\alpha},B}$. Denote by $M$ the Hardy-Littlewood maximal function
and $M^{k}$ the $k$-times iterations of $M$, then $M_{L(\log
L)^{k}}\approx M^{k+1}$ for $k=0,1,2,\cdots$. We also need the
following estimates in the proof of Theorem \ref{thm.1}.

\begin{lemma}\upcite{Z1}\label{lem.4}
Let $1\le p<\infty, \omega^{p}\in A_{1}$ and $B$ be a ball. Then for
any $y\in B$ and any positive integer $m$, there has
$$ \Big( \frac{1}{|2^{k}B|}  \dint_{2^{k}B} |b(x)-b_{B}|^{mp}\omega^{p}(x) \mathrm{d}x \Big)^{1/p}\le C \|b\|_{\ast}^{m} (k+1)^{m} \inf_{y\in B} \omega(y), \ k=0,1,2,\cdots.$$
\end{lemma}

\begin{lemma}\label{lem.5}
Let $1\le p<\infty, \omega^{p}\in A_{1}$ and $B$ be a ball. Then for
any $y\in B$ and any positive integer $m$, there has
$$ \Big( \frac{1}{|2^{k}B|}  \dint_{2^{k}B} \omega^{p}(x) \prod_{j=1}^{m} \Big|b_{j}(x)-(b_{j})_{B}\Big|^{p} \mathrm{d}x \Big)^{1/p}\le C \|\vec{b}\|_{\ast} (k+1)^{m} \inf_{y\in B} \omega(y), \ k=0,1,2,\cdots.$$
\end{lemma}

\begin{proof} By the H\"{o}lder's
inequality and Lemma \ref{lem.4}, we obtain
\begin{eqnarray*}
\Big( \frac{1}{|2^{k}B|}  \dint_{2^{k}B} \omega^{p}(x)
\prod_{j=1}^{m} \Big|b_{j}(x)-(b_{j})_{B}\Big|^{p} \mathrm{d}x
\Big)^{1/p} &\le & \prod_{j=1}^{m}  \Big( \frac{1}{|2^{k}B|}
\dint_{2^{k}B} \omega^{p}(x)
\Big|b_{j}(x)-(b_{j})_{B}\Big|^{p\gamma_{j}} \mathrm{d}x
\Big)^{\frac{1}{p\gamma_{j}}}\\
&\le & C\prod_{j=1}^{m}  \Big( \|b_{j}\|_{\ast}^{\gamma_{j}}
(k+1)^{\gamma_{j}} \inf_{y\in B} \omega(y)
\Big)^{\frac{1}{\gamma_{j}}}\\
&\le & C \|\vec{b}\|_{\ast} (k+1)^{m} \inf_{y\in B} \omega(y),
 \end{eqnarray*}
where $1=\frac{1}{\gamma_{1}} + \frac{1}{\gamma_{2}}+\cdots
+\frac{1}{\gamma_{m}}$.

 This completes the proof of Lemma \ref{lem.5}.
\end{proof}

We also need the following notations. For $\omega \in A_{_{\infty}}$
and a ball $B$, denote by
\begin{equation*}
\|f\|_{L(\log L)^{m},B,\omega}=\inf \Big\{ \lambda>0:
\frac{1}{\omega(B)}\dint_{B} \Phi_{m} \Big(\frac{|f(y)|}{\lambda}
\Big)\omega(y) \mathrm{d}y \le 1\Big\}
 \end{equation*}
and
\begin{equation*}
\|f\|_{\exp L^{1/m},B,\omega}=\inf \Big\{ \lambda>0:
\frac{1}{\omega(B)}\dint_{B} \tilde{\Phi}_{m}
\Big(\frac{|f(y)|}{\lambda} \Big)\omega(y) \mathrm{d}y \le 1\Big\}.
 \end{equation*}

Similar to (\ref{equ:3}), we have (c.f.\cite{RR}, p.69)
\begin{equation} \label{equ:7}
\|f\|_{L(\log L)^{m},B,\omega}\approx\inf \Big\{ \eta+
\frac{\eta}{\omega(B)}\dint_{B} \Phi_{m} \Big(\frac{|f(y)|}{\eta}
\Big)\omega(y) \mathrm{d}y\Big\}.
 \end{equation}

By (\ref{equ:4}), there also holds the following generalized
H\"{o}lder's inequality
\begin{equation} \label{equ:8}
\frac{1}{\omega(B)} \dint_{B} |f_{1}(y)\cdots f_{m}(y)g(y)|
\omega(y) \mathrm{d}y \le C\|g\|_{L(\log
L)^{m},B,\omega}\prod_{j=1}^{m}\|f_{j}\|_{\exp L,B,\omega}.
 \end{equation}

Furthermore, for any $b\in \BMO(\R^{n})$, any ball $B$ and any
$\omega \in A_{_{\infty}}$, there has
\begin{equation} \label{equ:9}
\|b-b_{B}\|_{\exp
L,B,\omega}\le C\|b\|_{\ast},
 \end{equation}

Indeed, by John-Nirenberg's inequality, there exist positive
constants $C_{1}$ and $C_{2}$, such that
$$|\{x\in B: |b(x)-b_{B}|>t\}|\le C_{1}|B| e^{-C_{2}t/\|b\|_{\ast}}.$$
Noting that $\omega \in A_{_{\infty}}$, from the proof of Theorem 5
in \cite{MW}, there is a $\delta>0$, such that
$$\omega(\{x\in B: |b(x)-b_{B}|>t\})\le C_{1}\omega(B) e^{-C_{2}\delta t/\|b\|_{\ast}}.$$

Similar to the proof of Corollary 7.1.7 in \cite{G} (p.528), we have
\begin{equation} \label{equ:10}
 \frac{1}{\omega(B)} \dint_{B} \exp \Big(\frac{|b(x)-b_{B}|}{C_{3}\|b\|_{\ast} } \Big)\omega(x)\mathrm{d}x \le
 C,
 \end{equation}
which implies (\ref{equ:9}).

\section{Proof of Theorem\ref{thm.1}}
\setcounter{equation}{0}

\begin{proof}  Without loss of generality, we may assume that for $j
=1,\ldots,m , \|b_{j}\|_{\Osc_{\exp L^{\tau_{_{j}}}}(\R^n)}=1$. In
fact, let $$ \tilde{b}_{j}=\frac{b_{j}}{\|b_{j}\|_{\Osc_{\exp L^{\tau_{_{j}}}}(\R^n)}}$$ for $j =1,\ldots,m$. The homogeneity
tells us that for any $\lambda>0$,
\begin{equation} \label{equ:3.1}
 \omega(\{x\in \R^{n}: \mu_{\Omega,\vec{b}}(f)(x)> \lambda\}) = \omega(\{x\in \R^{n}: \mu_{\Omega,\tilde{b}}(f)(x)> \lambda/\|\vec{b}\|_{\Osc_{\exp
L^{\tau}}(\R^n)}\})
 \end{equation}
Noting that $\|\tilde{b}_{j}\|_{\Osc_{\exp
L^{\tau_{_{j}}}}(\R^n)}=1$  for $j =1,\ldots,m$, if when
$\|b_{j}\|_{\Osc_{\exp L^{\tau_{_{j}}}}(\R^n)}=1~(j =1,\ldots,m)$,
the theorem is true. , by (\ref{equ:3.1}) and the inequality
$$ \Phi_{s}(t_{1} t_{2})\leq C \Phi_{s}(t_{1}) \Phi_{s}( t_{2})$$
for any $s>0, t_{1}, t_{2}\geq 0$, we easily obtain that the theorem
still holds for any $b_{j}\in \Osc_{\exp L^{\tau_{_{j}}}}(\R^n)~(j
=1,\ldots,m)$.

 For a fixed
$\lambda$, we consider the Calder\'{o}n-Zygmund decomposition of $f$
at height $\lambda$ and get a sequence of balls $\{B_{i}\}$, where
$B_{i}$ is a ball centered at $x_{i}$ with radius $r_{i}$, such that
$|f(x)|\le C\lambda$ for $ a.e.\ x\in \R^{n}\setminus\cup_{i}B_{i}$
and
\begin{equation} \label{equ:11}
 \lambda<  \frac{1}{|B_{i}|}\dint_{B_{i}} |f(y)| \mathrm{d}y\le 2^{n}\lambda.
 \end{equation}

Moreover, there is an integer $N\ge 1$, independent of $f$ and
$\lambda$, such that for every point in $\R^{n}$ belongs to at most
$N$ balls in $\{B_{i}\}$.

We decompose $f=g+h$, where
\begin{displaymath}
g(x)= \left\{\begin{array}{ll}
f(x), &  \ x\in \R^{n}\setminus\cup_{i}B_{i}, \\
f_{B_{i}}, & \  x\in B_{i}.
\end{array}\right.
\end{displaymath}
Then $h(x)= f(x)-g(x)=\sum_{i}h_{i}(x)$ with $h_{i}(x)=(f(x)-
f_{B_{i}}) \chi_{_{B_{i}}}(x)$. Obviously, $\supp h_{i}\subset
B_{i}, \int_{B_{i}} h_{i}(y) \mathrm{d}y=0$ and
\begin{equation} \label{equ:12}
 |g(x)|\le 2^{n}\lambda, \ \ \ \  a.e.\ x\in
\R^{n}.
 \end{equation}

Noting that if $\omega^{r'}\in A_{1}$ then $\omega\in A_{1}$, and
then $M(\omega)(x)\le C\omega(x)$ for $a.e.\ x\in \R^{n}$. By
(\ref{equ:11}) and the fact that $|B_{i}|^{-1}
\omega(B_{i})=|B_{i}|^{-1} \int_{B_{i}} \omega(x) \mathrm{d}x \le C
\inf_{y\in B_{i}} \omega(y)$, we have
\begin{equation} \label{equ:13}
\omega(B_{i}) \le C |B_{i}|\inf_{y\in B_{i}} \omega(y) \le
C\lambda^{-1}\int_{B_{i}} |f(y)| \mathrm{d}y \inf_{y\in B_{i}}
\omega(y)\le C\lambda^{-1}\int_{B_{i}} |f(y)| \omega(y) \mathrm{d}y
.
 \end{equation}

Denote by $E=\cup_{i}(4B_{i})$, it follows from (\ref{equ:13}) that
$$ \omega(E) \le C \sum_{i}\int_{B_{i}} \omega(x) \mathrm{d}x =C \sum_{i}\omega(B_{i}) \le C\lambda^{-1}\|f\|_{L^{1}_{\omega}(\R^{n})}.
$$

Write
 \begin{eqnarray*}
&&\;  \omega(\{x\in \R^{n}: \mu_{\Omega,\vec{b}}(f)(x)> \lambda\}) \le \omega(\{x\in \R^{n}\setminus E: \mu_{\Omega,\vec{b}}(f)(x)>  \lambda\})+ \omega(E)\\
&&\; \le \omega(\{x\in \R^{n}\setminus E:
\mu_{\Omega,\vec{b}}(g)(x)> \frac{\lambda}{2}\})
  + \omega(\{x\in \R^{n}\setminus E: \mu_{\Omega,\vec{b}}(h)(x)> \frac{\lambda}{2}\})+ \omega(E)\\
&&\; \le I_{1}+I_{2}+C\lambda^{-1}\|f\|_{L^{1}_{\omega}(\R^{n})}.
\end{eqnarray*}

We consider $I_{1}$ first. For $\omega^{r'}\in A_{1}$ there has
$\omega\in A_{1}$. Noting that $A_{1}\subset A_{s}~(s\ge 1)$, then
for any $p>r'$, we have $\omega\in A_{p/r'}$. It follows from Lemma
\ref{lem.3},  (\ref{equ:12}) and (\ref{equ:13}) that
\begin{equation} \label{equ:14}
\begin{split}
I_{1} & 
\le C \lambda^{-p}  \dint_{\R^{n}} \Big( \mu_{\Omega,\vec{b}}(g)(x) \Big)^{p} \omega(x) \mathrm{d}x   
\le C 
       \lambda^{-p}  \dint_{\R^{n}} |g(x)|^{p} \omega(x) \mathrm{d}x  \le C  
                                                           \lambda^{-1}  \dint_{\R^{n}} |g(x)| \omega(x) \mathrm{d}x \\
&\le C 
          \lambda^{-1}  \Big( \dint_{\R^{n}\setminus\cup_{i}B_{i}} |g(x)| \omega(x) \mathrm{d}x + \dint_{\cup_{i}B_{i}} |g(x)| \omega(x) \mathrm{d}x \Big)\\
 &\le C 
               \lambda^{-1}  \Big( \dint_{\R^{n}} |f(x)| \omega(x) \mathrm{d}x + \sum_{i}\dint_{B_{i}} |f_{B_{i}}| \omega(x) \mathrm{d}x \Big)\\
 &\le C 
                  \lambda^{-1}\|f\|_{L^{1}_{\omega}(\R^{n})}+ C 
                                         \lambda^{-1}   \sum_{i}\dint_{B_{i}} \Big(|B_{i}|^{-1}\dint_{B_{i}}|f(y)|\mathrm{d}y \Big)  \omega(x) \mathrm{d}x \\
&\le C 
        \lambda^{-1}\|f\|_{L^{1}_{\omega}(\R^{n})}+ C 
                                     \lambda^{-1}   \sum_{i}\dint_{B_{i}}|f(y)|\mathrm{d}y \Big(|B_{i}|^{-1}\dint_{B_{i}}  \omega(x) \mathrm{d}x \Big) \\
&\le C 
             \lambda^{-1}\|f\|_{L^{1}_{\omega}(\R^{n})}+ C 
                                               \lambda^{-1}   \sum_{i}\dint_{B_{i}}|f(y)|\mathrm{d}y  \inf_{y\in B_{i}}  \omega(y) \\
&\le C  
                \lambda^{-1}\|f\|_{L^{1}_{\omega}(\R^{n})}+ C 
                                         \lambda^{-1}   \sum_{i}\dint_{B_{i}}|f(y)| \omega(y) \mathrm{d}y   \\
&\le C 
 \lambda^{-1}\|f\|_{L^{1}_{\omega}(\R^{n})}.
 \end{split}
 \end{equation}

We remark that the proof of (\ref{equ:14}) implies the following
fact, which will be used later.
\begin{equation} \label{equ:15}
\sum_{i}\dint_{B_{i}} |f_{B_{i}}| \omega(x) \mathrm{d}x \le
C\|f\|_{L^{1}_{\omega}(\R^{n})}.
 \end{equation}

Now, let us estimate $I_{2}$. By the definition of $\mu_{\Omega}$
and $\mu_{\Omega,\vec{b}}$, with the aid of the formula
\begin{eqnarray*}
 \prod_{j=1}^{m} \big(  b_{j}(x)-b_{j}(y) \big) = \sum_{j=0}^{m}\sum_{\sigma\in \mathscr{C}_{j}^{m}} \Big(b(x)-b_{B_{i}} \Big)_{\sigma}  \Big(b_{B_{i}} -b(y) \Big)_{\sigma'},
             \end{eqnarray*}
we have
\begin{eqnarray*}
&&\;\mu_{\Omega,\vec{b}}(h)(x) 
= \Big(\dint_{0}^{\infty} \Big| \dint_{|x-y|\le t} \frac{\Omega(x-y)h(y) }{|x-y|^{n-1}} \sum_{j=0}^{m}\sum_{\sigma\in \mathscr{C}_{j}^{m}}  \Big(b(x)-b_{B_{i}} \Big)_{\sigma}  \Big(b_{B_{i}} -b(y) \Big)_{\sigma'}\mathrm{d}y\Big|^{2} \frac{\mathrm{d}t}{t^{3}}\Big)^{\frac{1}{2}} \\
 &&\;\le  \Big(\dint_{0}^{\infty} \Big| \dint_{|x-y|\le t} \frac{\Omega(x-y)h(y) }{|x-y|^{n-1}} \prod_{j=1}^{m}  \Big(b_{j}(x)-(b_{j})_{B_{i}} \Big)\mathrm{d}y\Big|^{2} \frac{\mathrm{d}t}{t^{3}}\Big)^{\frac{1}{2}} \\
 &&\; \ \ \ + \Big(\dint_{0}^{\infty} \Big| \dint_{|x-y|\le t} \frac{\Omega(x-y)h(y) }{|x-y|^{n-1}} \sum_{j=1}^{m-1}\sum_{\sigma\in \mathscr{C}_{j}^{m}}  \Big(b(x)-b_{B_{i}} \Big)_{\sigma}  \Big(b_{B_{i}} -b(y) \Big)_{\sigma'}\mathrm{d}y\Big|^{2} \frac{\mathrm{d}t}{t^{3}}\Big)^{\frac{1}{2}} \\
 &&\;\ \ \ +\Big(\dint_{0}^{\infty} \Big| \dint_{|x-y|\le t} \frac{\Omega(x-y)h(y) }{|x-y|^{n-1}} \prod_{j=1}^{m}  \Big((b_{j})_{B_{i}} -b_{j}(y) \Big)\mathrm{d}y\Big|^{2} \frac{\mathrm{d}t}{t^{3}}\Big)^{\frac{1}{2}} \\
 &&\;\le  \sum_{i}\prod_{j=1}^{m}  \Big|b_{j}(x)-(b_{j})_{B_{i}} \Big| \mu_{\Omega}(h_{i})(x)   + \sum_{i} \sum_{j=1}^{m-1}\sum_{\sigma\in \mathscr{C}_{j}^{m}}  \Big|\big(b(x)-b_{B_{i}} \big)_{\sigma} \Big|  \mu_{\Omega} \big(h_{i}  \big(b_{B_{i}} -b \big)_{\sigma'} \big)(x) \\
 &&\;\ \ \ +\mu_{\Omega} \big(\sum_{i} h_{i} \prod_{j=1}^{m}  \big((b_{j})_{B_{i}} -b_{j} \big)
 \big)(x).
\end{eqnarray*}
So, we can write $I_{2}$ as
\begin{equation} \label{equ:16}
\begin{split}
I_{2} 
 &\le  \omega\Big( \Big\{x\in \R^{n}\setminus E: \sum_{i}\prod_{j=1}^{m}  \Big|b_{j}(x)-(b_{j})_{B_{i}} \Big| \mu_{\Omega}(h_{i})(x) > \frac{\lambda}{6} \Big\} \Big) \\
& \ + \omega\Big( \Big\{x\in \R^{n}\setminus E: \sum_{i} \sum_{j=1}^{m-1}\sum_{\sigma\in \mathscr{C}_{j}^{m}}  \Big|\big(b(x)-b_{B_{i}} \big)_{\sigma} \Big|  \mu_{\Omega} \big(h_{i}  \big(b_{B_{i}} -b \big)_{\sigma'} \big)(x)> \frac{\lambda}{6}\Big\} \Big) \\
 &\ + \omega\Big( \Big\{x\in \R^{n}\setminus E: \mu_{\Omega} \big(\sum_{i}h_{i} \prod_{j=1}^{m}  \big((b_{j})_{B_{i}} -b_{j} \big)
 \big)(x)> \frac{\lambda}{6} \Big\} \Big) \\
 &=I_{21}+I_{22}+I_{23}.
\end{split}
 \end{equation}

For $I_{21}$, using chebychev's inequality and Minkowski's
inequality, we have
\begin{equation} \label{equ:17}
\begin{split}
&I_{21}= \omega\Big( \Big\{x\in \R^{n}\setminus E: \sum_{i}\prod_{j=1}^{m}  \Big|b_{j}(x)-(b_{j})_{B_{i}} \Big| \mu_{\Omega}(h_{i})(x) > \frac{\lambda}{6} \Big\} \Big) \\
 &\le  C\lambda^{-1} \sum_{i} \dint_{\R^{n}\setminus 4B_{i}}\prod_{j=1}^{m}  \Big|b_{j}(x)-(b_{j})_{B_{i}} \Big| \mu_{\Omega}(h_{i})(x) \omega(x) \mathrm{d}x \\
 &\le  C\lambda^{-1} \sum_{i} \dint_{\R^{n}\setminus 4B_{i}}\prod_{j=1}^{m}  \Big|b_{j}(x)-(b_{j})_{B_{i}} \Big|  \Big(\dint_{0}^{|x-x_{i}|+2r_{i}} \Big| \dint_{|x-y|\le t} \frac{\Omega(x-y)h_{i}(y) }{|x-y|^{n-1}} \mathrm{d}y\Big|^{2} \frac{\mathrm{d}t}{t^{3}}\Big)^{\frac{1}{2}} \omega(x) \mathrm{d}x \\
& \ +  C\lambda^{-1} \sum_{i} \dint_{\R^{n}\setminus 4B_{i}}\prod_{j=1}^{m}  \Big|b_{j}(x)-(b_{j})_{B_{i}} \Big|  \Big(\dint_{|x-x_{i}|+2r_{i}}^{\infty} \Big| \dint_{|x-y|\le t} \frac{\Omega(x-y)h_{i}(y) }{|x-y|^{n-1}} \mathrm{d}y\Big|^{2} \frac{\mathrm{d}t}{t^{3}}\Big)^{\frac{1}{2}} \omega(x) \mathrm{d}x \\
 &=I_{211}+I_{212}.
\end{split}
 \end{equation}

For $x\in \R^{n}\setminus 4B_{i}$ and $y\in B_{i}$, there has
$|x-y|\le |x-x_{i}|+r_{i}$ and $|x-y|\sim |x-x_{i}|\sim
|x-x_{i}|+2r_{i}$, and then
$$ \dint_{|x-y|}^{|x-x_{i}|+2r_{i}} \frac{\mathrm{d}t}{t^{3}}=\frac{1}{2}\Big(\frac{1}{|x-y|^{2}} - \frac{1}{(|x-x_{i}|+2r_{i})^{2}}\Big) \le \frac{Cr_{i}}{|x-y|^{3}}. $$

Noting that $\supp h_{i} \subset B_{i}$, it follows from the
Minkowski's inequality that
\begin{equation*}
\begin{split}
&I_{211} 
 \le  C\lambda^{-1} \sum_{i} \dint_{\R^{n}\setminus 4B_{i}}\prod_{j=1}^{m}  \Big|b_{j}(x)-(b_{j})_{B_{i}} \Big|  \Big(\dint_{B_{i}}  \frac{|\Omega(x-y)||h_{i}(y)| }{|x-y|^{n-1}} \Big(\dint_{|x-y|}^{|x-x_{i}|+2r_{i}}  \frac{\mathrm{d}t}{t^{3}}\Big)^{\frac{1}{2}} \mathrm{d}y \Big)\omega(x) \mathrm{d}x \\
\end{split}
 \end{equation*}
 \begin{equation}\label{equ:18}
\begin{split} &\le   C\lambda^{-1} \sum_{i} r_{i}^{1/2}\dint_{\R^{n}\setminus 4B_{i}}\prod_{j=1}^{m}  \Big|b_{j}(x)-(b_{j})_{B_{i}} \Big|  \Big(\dint_{B_{i}}  \frac{|\Omega(x-y)||h_{i}(y)| }{|x-y|^{n+1/2}}  \mathrm{d}y \Big)\omega(x) \mathrm{d}x \\
&\le   C\lambda^{-1} \sum_{i} r_{i}^{1/2}\dint_{B_{i}}|h_{i}(y)| \sum_{k=1}^{\infty} \Big(\dint_{2^{k+1}B_{i}\setminus 2^{k}B_{i}}  \frac{|\Omega(x-y)|}{|x-y|^{n+1/2}}  \prod_{j=1}^{m}  \Big|b_{j}(x)-(b_{j})_{B_{i}} \Big| \omega(x)\mathrm{d}x \Big) \mathrm{d}y \\
 & \le   C\lambda^{-1} \sum_{i} r_{i}^{1/2}\dint_{B_{i}}|h_{i}(y)| \sum_{k=1}^{\infty}\bigg( \Big(\dint_{2^{k+1}B_{i}\setminus 2^{k}B_{i}}  \frac{|\Omega(x-y)|^{r}}{|x-y|^{n+1/2}} \mathrm{d}x \Big)^{1/r}\\
 &  \ \ \ \hspace{3cm}  \times \Big(\dint_{2^{k+1}B_{i}\setminus 2^{k}B_{i}}\frac{  \omega^{r'}(x)}{|x-y|^{n+1/2}}\prod_{j=1}^{m}  \Big|b_{j}(x)-(b_{j})_{B_{i}} \Big|^{r'}\mathrm{d}x \Big)^{1/r'} \bigg) \mathrm{d}y \\
\end{split}
 \end{equation}

Noting that $2^{k-1}r_{i}\le |x-y|\le 2^{k+2}r_{i}$ whenever $y\in
B_{i}$ and $x\in 2^{k+1}B_{i}\setminus 2^{k}B_{i}$, we have
\begin{equation} \label{equ:19}
\begin{split}
&  \Big(\dint_{2^{k+1}B_{i}\setminus 2^{k}B_{i}}
\frac{|\Omega(x-y)|^{r}}{|x-y|^{n+1/2}} \mathrm{d}x \Big)^{1/r}
 \le   \Big(\dint_{2^{k-1}r_{i}\le |x-y|\le 2^{k+2}r_{i}}  \frac{|\Omega(x-y)|^{r}}{|x-y|^{n+1/2}} \mathrm{d}x \Big)^{1/r} \\
& \le   \Big(\dint_{2^{k-1}r_{i} }^{2^{k+2}r_{i}} \rho^{n-1}
\Big(\dint_{S^{n-1}}\frac{|\Omega(x')|^{r}}{\rho^{n+1/2}}
\mathrm{d}x' \Big)\mathrm{d}\rho \Big)^{1/r}
 \le   C (2^{k}r_{i})^{-\frac{1}{2r}} \|\Omega\|_{L^{r}(S^{n-1})}.
 \end{split}
 \end{equation}

And noting that $\omega^{r'}\in A_{1}$ and $\|b_{j}\|_{\BMO}\le C
\|b_{j}\|_{\Osc_{\exp L^{\tau_{_{j}}}}}$ for $\tau_{_{j}}\ge 1
~(1\le j\le m)$, by the H\"{o}lder's inequality, Minkowski's
inequality, the properties of $\BMO(\R^{n})$ functions and Lemma
\ref{lem.5},
 we have
\begin{equation} \label{equ:20}
\begin{split}
&\Big(\dint_{2^{k+1}B_{i}\setminus 2^{k}B_{i}}\frac{ \omega^{r'}(x)}{|x-y|^{n+1/2}}\prod_{j=1}^{m} \Big|b_{j}(x)-(b_{j})_{B_{i}} \Big|^{r'}\mathrm{d}x \Big)^{1/r'} \\
 &\le C (2^{k+1}r_{i})^{-(n+\frac{1}{2})/r'}\Big(\dint_{2^{k+1}B_{i}} \omega^{r'}(x)\prod_{j=1}^{m} \Big|b_{j}(x)-(b_{j})_{B_{i}} \Big|^{r'}\mathrm{d}x \Big)^{1/r'}  \\
 &\le C (2^{k+1}r_{i})^{-(n+\frac{1}{2})/r'}\Big( \frac{|2^{k+1}B_{i}|}{|2^{k+1}B_{i}|}\dint_{2^{k+1}B_{i}} \omega^{r'}(x)\prod_{j=1}^{m} \Big|b_{j}(x)-(b_{j})_{B_{i}} \Big|^{r'}\mathrm{d}x \Big)^{1/r'}  \\
&\le   C 
  (2^{k}r_{i})^{-\frac{1}{2r'}} (k+1)^{m} \inf_{y\in B_{i}}  \omega(y) .\\
\end{split}
 \end{equation}

This, together with  (\ref{equ:18}) and (\ref{equ:19}), gives
\begin{equation} \label{equ:21}
\begin{split}
&I_{211}\le   C \|\Omega\|_{L^{r}(S^{n-1})} 
                 \lambda^{-1} \sum_{i} r_{i}^{1/2}\dint_{B_{i}}|h_{i}(y)| \Big(\sum_{k=1}^{\infty}(k+1)^{m} (2^{k}r_{i})^{-\frac{1}{2}}\Big)\omega(y) \mathrm{d}y \\
 &\le   C 
                          \lambda^{-1} \sum_{i} \dint_{B_{i}}|h_{i}(y)| \Big(\sum_{k=1}^{\infty}(k+1)^{m} 2^{-k/2}\Big)\omega(y) \mathrm{d}y
 \le   C 
                      \lambda^{-1} \sum_{i} \dint_{B_{i}}|h_{i}(y)| \omega(y) \mathrm{d}y. \\
 \end{split}
 \end{equation}

Next, let us consider $I_{212}$. Write
$K(x,y,x_{i})=\frac{\Omega(x-y)}{|x-y|^{n-1}}-\frac{\Omega(x-x_{i})}{|x-x_{i}|^{n-1}}$
for simplicity. Noting that for any $y\in B_{i}$, any $x\in
\R^{n}\setminus 4 B_{i}$ and $t$ with $|x-x_{i}|+2r_{i}\le t$, there
has $|x-y|\le |x-x_{i}|+r_{i}< t$, then by the cancellation
condition of $h_{i}$, we have
\begin{eqnarray*}
&&\;I_{212} 
 \le  C\lambda^{-1} \sum_{i} \dint_{\R^{n}\setminus 4B_{i}}\prod_{j=1}^{m}  \Big|b_{j}(x)-(b_{j})_{B_{i}} \Big|  \Big(\dint_{B_{i}}  |K(x,y,x_{i})||h_{i}(y)| \Big(\dint_{|x-x_{i}|+2r_{i}}^{\infty}  \frac{\mathrm{d}t}{t^{3}}\Big)^{\frac{1}{2}} \mathrm{d}y \Big)\omega(x) \mathrm{d}x \\
  &&\;\le   C\lambda^{-1} \sum_{i}\dint_{\R^{n}\setminus 4B_{i}}\prod_{j=1}^{m}  \Big|b_{j}(x)-(b_{j})_{B_{i}} \Big|  \Big(\dint_{B_{i}}  \frac{|K(x,y,x_{i})||h_{i}(y)| }{|x-x_{i}|}  \mathrm{d}y \Big)\omega(x) \mathrm{d}x \\
&&\;\le   C\lambda^{-1} \sum_{i} \dint_{B_{i}}|h_{i}(y)| \sum_{k=1}^{\infty} (2^{k}r_{i})^{-1}\Big(\dint_{2^{k+1}B_{i}\setminus 2^{k}B_{i}}  |K(x,y,x_{i})|  \prod_{j=1}^{m}  \Big|b_{j}(x)-(b_{j})_{B_{i}} \Big| \omega(x)\mathrm{d}x \Big) \mathrm{d}y \\
\end{eqnarray*}

By the H\"{o}lder's inequality, Lemma \ref{lem.1} and Lemma
\ref{lem.5}, there has
\begin{eqnarray*}
&&\;\dint_{2^{k+1}B_{i}\setminus 2^{k}B_{i}}  |K(x,y,x_{i})|  \prod_{j=1}^{m}  \Big|b_{j}(x)-(b_{j})_{B_{i}} \Big| \omega(x)\mathrm{d}x \\
 &&\;\le \Big( \dint_{2^{k+1}B_{i}\setminus 2^{k}B_{i}}  |K(x,y,x_{i})|^{r} \mathrm{d}x \Big)^{1/r} \Big( \dint_{2^{k+1}B_{i}\setminus 2^{k}B_{i}} \prod_{j=1}^{m}  \Big|b_{j}(x)-(b_{j})_{B_{i}} \Big|^{r'} \omega^{r'}(x)\mathrm{d}x \Big)^{1/r'}\\
 &&\;\le   C 
              (k+1)^{m}2^{k}r_{i}\Big(2^{-k} + \dint_{\frac{|y-x_{i}|}{2^{k+1}r_{i}}}^{\frac{|y-x_{i}|}{2^{k}r_{i}}}\frac{\omega_{r}(\delta)}{\delta} \mathrm{d}\delta \Big) \inf_{y\in B_{i}}\omega(y). \\
\end{eqnarray*}

Therefore,
\begin{equation} \label{equ:22}
\begin{split}
&I_{212}\le   C 
              \lambda^{-1} \sum_{i} \dint_{B_{i}}|h_{i}(y)| \omega(y)\sum_{k=1}^{\infty}  (k+1)^{m}\Big(2^{-k} + \dint_{\frac{|y-x_{i}|}{2^{k+1}r_{i}}}^{\frac{|y-x_{i}|}{2^{k}r_{i}}}\frac{\omega_{r}(\delta)}{\delta} \mathrm{d}\delta \Big) \mathrm{d}y \\
 &\le   C 
            \lambda^{-1} \sum_{i} \dint_{B_{i}}|h_{i}(y)| \omega(y)\Big(\sum_{k=1}^{\infty}  (k+1)^{m}2^{-k} + \dint_{0}^{1}\frac{\omega_{r}(\delta)}{\delta}\big(\log\frac{1}{\delta} \big)^{m} \mathrm{d}\delta \Big) \mathrm{d}y \\
 &\le   C 
              \lambda^{-1} \sum_{i} \dint_{B_{i}}|h_{i}(y)| \omega(y) \mathrm{d}y. \\
 \end{split}
 \end{equation}

Note that $h_{i}(y)=f(y)+f_{B_{i}}$ when $y\in B_{i}$, it follows
from (\ref{equ:15}), (\ref{equ:17}), (\ref{equ:21}) and
(\ref{equ:22}) that
\begin{eqnarray*}
I_{21}&\le &C 
\lambda^{-1} \sum_{i} \dint_{B_{i}}|h_{i}(y)| \omega(y) \mathrm{d}y
\le C 
             \lambda^{-1} \sum_{i} \dint_{B_{i}}( |f(y)|+|f_{B_{i}}|) \omega(y) \mathrm{d}y  \le C 
             \lambda^{-1}\|f\|_{L^{1}_{\omega}(\R^{n})}.
\end{eqnarray*}

To estimate $I_{23}$, noting that $\Omega\in L^{r}(S^{n-1})$ for
some $r>1$ and $\omega^{r'}\in A_{1}$,   using Lemma \ref{lem.2},
(\ref{equ:8}), (\ref{equ:9}), Lemma \ref{lem.5}, (\ref{equ:7}) and
(\ref{equ:13}), we have
\begin{eqnarray*}
I_{23} 
&\le&  \omega\Big( \Big\{x\in \R^{n}: \mu_{\Omega} \big(\sum_{i}h_{i} \prod_{j=1}^{m}  \big((b_{j})_{B_{i}} -b_{j} \big)  \big)(x)> \frac{\lambda}{6} \Big\} \Big) \\
&\le& C \lambda^{-1} \dint_{\R^{n}}  \sum_{i} |h_{i}(x)| \omega(x) \prod_{j=1}^{m}  \big|(b_{j})_{B_{i}} -b_{j}(x) \big|\mathrm{d}x \\
&\le& C \lambda^{-1}\sum_{i} \Big( \dint_{B_{i}}  |f(x)| \omega(x) \prod_{j=1}^{m}  \big|(b_{j})_{B_{i}} -b_{j}(x) \big|\mathrm{d}x +  \dint_{B_{i}}  |f_{B_{i}}| \omega(x) \prod_{j=1}^{m}  \big|(b_{j})_{B_{i}} -b_{j}(x) \big|\mathrm{d}x \Big)\\
&\le& C \lambda^{-1}\sum_{i} \omega(B_{i}) \|f\|_{L(\log L)^{m},B_{i},\omega}   \prod_{j=1}^{m} \big\|  |b_{j}-(b_{j})_{B_{i}} |\big\|_{\exp L,B_{i},\omega}  \\
&\ & \ \  +  C \lambda^{-1}\sum_{i}  \frac{1}{|B_{i}|} \dint_{B_{i}} |f(y)|\mathrm{d}y \dint_{B_{i}}  \omega(x) \prod_{j=1}^{m}  \big|(b_{j})_{B_{i}} -b_{j}(x) \big|\mathrm{d}x \\
\end{eqnarray*}
\begin{eqnarray*}&\le& C 
                         \lambda^{-1}\sum_{i} \Big( \omega(B_{i}) \|f\|_{L(\log L)^{m},B_{i},\omega} +   \dint_{B_{i}} |f(y)|\mathrm{d}y   \inf_{y\in B_{i}}\omega(y) \Big) \\
&\le& C 
         \lambda^{-1}\sum_{i} \Big( \omega(B_{i}) \inf \Big\{\lambda+\frac{\lambda}{\omega(B_{i})} \dint_{B_{i}} \Phi_{m}\Big(\frac{|f(y)|}{\lambda}
\Big)\omega(y) \mathrm{d}y \Big\}  +    \dint_{B_{i}} |f(y)|\omega(y)\mathrm{d}y  \Big)  \\
&\le& C 
        \sum_{i} \Big( \omega(B_{i})+\dint_{B_{i}} \Phi_{m}\Big(\frac{|f(y)|}{\lambda}
\Big)\omega(y) \mathrm{d}y \Big)  +    C 
                                      \lambda^{-1}\dint_{\R^{n}} |f(y)|\omega(y)\mathrm{d}y   \\
&\le& C 
         \sum_{i} \Big( \lambda^{-1}\int_{B_{i}} |f(y)| \omega(y) \mathrm{d}y+\dint_{B_{i}} \frac{|f(y)|}{\lambda}\Big(1+\log^{+}\frac{|f(y)|}{\lambda} \Big)^{m}\omega(y) \mathrm{d}y \Big)
+    C 
             \lambda^{-1}\dint_{\R^{n}} |f(y)|\omega(y)\mathrm{d}y   \\
&\le& C 
          \dint_{\R^{n}} \frac{|f(y)|}{\lambda}\Big(1+\log^{+}\frac{|f(y)|}{\lambda} \Big)^{m}\omega(y) \mathrm{d}y . \\
\end{eqnarray*}

 Now, let us turn to estimate for $I_{22}$. Using the Minkowski's
inequality,  we have
\begin{eqnarray*}
&&\;I_{22}=  \omega\Big( \Big\{x\in \R^{n}\setminus E: \sum_{i} \sum_{j=1}^{m-1}\sum_{\sigma\in \mathscr{C}_{j}^{m}}  \Big|\big(b(x)-b_{B_{i}} \big)_{\sigma} \Big|  \mu_{\Omega} \big(h_{i}  \big(b_{B_{i}} -b \big)_{\sigma'} \big)(x)> \frac{\lambda}{6}\Big\} \Big) \\
&&\;\le  C\lambda^{-1} \sum_{i} \sum_{j=1}^{m-1}\sum_{\sigma\in \mathscr{C}_{j}^{m}}  \dint_{\R^{n}\setminus 4B_{i}}\Big|\big(b(x)-b_{B_{i}} \big)_{\sigma} \Big|  \mu_{\Omega} \big(h_{i}  \big(b_{B_{i}} -b \big)_{\sigma'} \big)(x) \omega(x) \mathrm{d}x \\
&&\;\le  C\lambda^{-1} \sum_{i} \sum_{j=1}^{m-1}\sum_{\sigma\in \mathscr{C}_{j}^{m}}  \dint_{\R^{n}\setminus 4B_{i}}\Big|\big(b(x)-b_{B_{i}} \big)_{\sigma} \Big|  \Big(\dint_{0}^{|x-x_{i}|+2r_{i}} \Big| \dint_{|x-y|\le t} \frac{\Omega(x-y)h_{i}(y) }{|x-y|^{n-1}} \\
 &&\;\ \ \ \ \  \times    \big(b_{B_{i}} -b(y) \big)_{\sigma'}  \mathrm{d}y\Big|^{2} \frac{\mathrm{d}t}{t^{3}}\Big)^{\frac{1}{2}} \omega(x) \mathrm{d}x \\
&&\; \  +   C\lambda^{-1} \sum_{i} \sum_{j=1}^{m-1}\sum_{\sigma\in \mathscr{C}_{j}^{m}}  \dint_{\R^{n}\setminus 4B_{i}}\Big|\big(b(x)-b_{B_{i}} \big)_{\sigma} \Big|   \Big(\dint_{|x-x_{i}|+2r_{i}}^{\infty} \Big| \dint_{|x-y|\le t} \frac{\Omega(x-y)h_{i}(y) }{|x-y|^{n-1}}\\
    &&\;\ \ \ \ \  \times      \big(b_{B_{i}} -b(y) \big)_{\sigma'}  \mathrm{d}y\Big|^{2} \frac{\mathrm{d}t}{t^{3}}\Big)^{\frac{1}{2}} \omega(x) \mathrm{d}x = C \lambda^{-1} \sum_{i} 
               (I_{221}+ I_{222})  . \\
\end{eqnarray*}

For $I_{221}$ and $I_{222}$,  similar to the estimates for $I_{21}$
and $I_{23}$, we can get
\begin{eqnarray*}
&&\;I_{221} 
\le  C \sum_{j=1}^{m-1}\sum_{\sigma\in \mathscr{C}_{j}^{m}} r_{i}^{1/2}\dint_{\R^{n}\setminus 4B_{i}}\Big|\big(b(x)-b_{B_{i}} \big)_{\sigma} \Big|  \Big(\dint_{B_{i}}  \frac{|\Omega(x-y)| |h_{i}(y)| }{|x-y|^{n+1/2}}    \big(b_{B_{i}} -b(y) \big)_{\sigma'}  \mathrm{d}y\Big) \omega(x) \mathrm{d}x \\
&&\;\le C 
\Big(\omega(B_{i}) \inf
\Big\{\lambda+\frac{\lambda}{\omega(B_{i})} \dint_{B_{i}}
\Phi_{m}\Big(\frac{|f(y)|}{\lambda}
\Big)\omega(y) \mathrm{d}y \Big\}  +    \dint_{B_{i}} |f(y)|\omega(y)\mathrm{d}y  \Big). \\
\\
&&\;I_{222}
\le  C\sum_{j=1}^{m-1}\sum_{\sigma\in \mathscr{C}_{j}^{m}} \dint_{\R^{n}\setminus 4B_{i}}\Big|\big(b(x)-b_{B_{i}} \big)_{\sigma} \Big|  \Big(\dint_{B_{i}}  \frac{|K(x,y,x_{i})| |h_{i}(y)| }{|x-x_{i}|}    \big(b_{B_{i}} -b(y) \big)_{\sigma'}  \mathrm{d}y\Big) \omega(x) \mathrm{d}x \\
&&\;\le C 
  \Big(\omega(B_{i}) \inf \Big\{\lambda+\frac{\lambda}{\omega(B_{i})} \dint_{B_{i}}
\Phi_{m}\Big(\frac{|f(y)|}{\lambda}
\Big)\omega(y) \mathrm{d}y \Big\}  +    \dint_{B_{i}} |f(y)|\omega(y)\mathrm{d}y  \Big). \\
\end{eqnarray*}

Thus, we have
\begin{eqnarray*}
 &&\;I_{22}\le  C \lambda^{-1} \sum_{i} 
                        \Big(\omega(B_{i}) \inf \Big\{\lambda+\frac{\lambda}{\omega(B_{i})} \dint_{B_{i}} \Phi_{m}\Big(\frac{|f(y)|}{\lambda} \Big)\omega(y) \mathrm{d}y \Big\}
 +    \dint_{B_{i}} |f(y)|\omega(y)\mathrm{d}y  \Big) \\ 
&&\; \le   C 
             \dint_{\R^{n}} \frac{|f(y)|}{\lambda}\Big(1+\log^{+}\frac{|f(y)|}{\lambda}
\Big)^{m}\omega(y) \mathrm{d}y .
\end{eqnarray*}

From (\ref{equ:16}) and the above estimates for $I_{21}, I_{22}$ and
$I_{23}$, we have
$$ I_{2} \le C 
          \dint_{\R^{n}}
\frac{|f(y)|}{\lambda}\Big(1+\log^{+}\frac{|f(y)|}{\lambda}
\Big)^{m}\omega(y) \mathrm{d}y .$$

This finishes the proof of Theorem \ref{thm.1}.
\end{proof}


\bigskip

\end{document}